\documentclass[authoryear, preprint,11pt]{elsarticle}

\usepackage{amssymb}
\usepackage{url}
\usepackage{graphicx}
\usepackage{amsmath}
\usepackage[gen]{eurosym}
\usepackage{wrapfig}
\usepackage{tabularx}
\usepackage{booktabs}
\usepackage{tikz}
\tikzset{
	treenode/.style = {shape=rectangle, rounded corners,
		draw, align=center,
		top color=white, bottom color=blue!20},
	root/.style     = {treenode, font=\Large, bottom color=red!30},
	env/.style      = {treenode, font=\ttfamily\normalsize},
	dummy/.style    = {circle,draw}
}
\usepackage{amssymb}
\usepackage{subfig}

\usepackage{enumitem}
\usepackage[gen]{eurosym}
\usetikzlibrary{plotmarks}
\usepackage{adjustbox}
\usepackage{ragged2e}

\newtheorem{thm}{Theorem}
\newtheorem{defi}[thm]{Definition}
\newtheorem{lem}[thm]{Lemma}

\newtheorem{rem}[thm]{Remark}

\newproof{pf}{Proof}

\DeclareMathOperator*{\argmax}{arg\,max}
\DeclareMathOperator*{\argmin}{arg\,min}

\usepackage{lipsum}
\makeatletter
\def\ps@pprintTitle{%
 \let\@oddhead\@empty
 \let\@evenhead\@empty
 \def\@oddfoot{}%
 \let\@evenfoot\@oddfoot}
\makeatother
\begin{document}

\begin{frontmatter}



\title{Multicriteria Adjustable Robustness}


\author[label1]{Elisabeth Halser\corref{cor1}}
\ead{elisabeth.halser@itwm.fraunhofer.de}
\author[label1]{Elisabeth Finhold}
\ead{elisabeth.finhold@itwm.fraunhofer.de}
\author[label1]{Neele Leith\"auser}
\ead{neele.leithaeuser@itwm.fraunhofer.de}
\author[label1]{Jan Schwientek}
\ead{jan.schwientek@itwm.fraunhofer.de}
\author[label1]{Katrin Teichert}
\ead{katrin.teichert@itwm.fraunhofer.de}
\author[label1]{Karl\nobreakdash-Heinz K\"ufer}
\ead{karl-heinz.kuefer@itwm.fraunhofer.de}

\cortext[cor1]{Corresponding Author}

\affiliation[label1]{organization={Fraunhofer Institute for Industrial Mathematics ITWM},
            addressline={Fraunhofer-Platz~1}, 
            city={Kaiserslautern},
            postcode={67663}, 
            country={Germany}}

\begin{abstract}
	Multicriteria adjustable robust optimization (MARO) problems arise in a wide variety of practical settings, for example, in the design of a building's energy supply. However, no general approaches, neither for the characterization of solutions to this problem class, nor potential solution methods, are available in the literature so far. We give different definitions for efficient solutions to MARO problems and look at three computational concepts to deal with the problems. These computational concepts can also be understood as additional solution definitions. We assess the advantages and disadvantages of the different computational approaches and analyze their connections to our initial definitions of MARO-efficiency. We observe that an $\varepsilon$-constraint inspired first-scalarize-then-robustify computational approach is beneficial because it provides an efficient set that is easy to understand for decision makers and provides tight bounds on the worst-case evaluation for a particular efficient solution. In contrast, a weighted sum first-scalarize-then-robustify approach keeps the problem structure more simple but is only beneficial if the desired trade-off between objectives is already known because the efficient set might look ambiguous. Further, we demonstrate that a first-robustify procedure only gives bad bounds and can be too optimistic as well as too pessimistic.
\end{abstract}

\begin{graphicalabstract}
\includegraphics{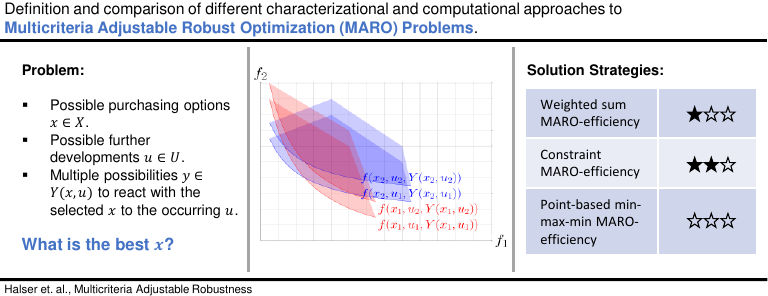}
\end{graphicalabstract}
%
\begin{highlights}
\item Characterizations of efficient solutions to multicriteria adjustable robust problems
\item Computational approaches to multicriteria adjustable robust problems
\item $\varepsilon$-constraint inspired approach has most interpretable bounds and intuitive solution
\item Weighted sum approach and directly-scalarized approach behave worse in general
\end{highlights}


\end{frontmatter}


\section{Introduction}

In practice, it is very common that purchasing decisions have to be made before an uncertain parameter reveals itself while afterwards it is still possible to react to this parameter in a certain way. This is the case, for example, for cost-optimal building energy supply design, where heating and cooling units have to be purchased without knowing the weather and price fluctuations in the future. However, it is still possible to react to these uncertainties when operating the units, once the weather and prices are known. These problems are known in the literature as adaptive or adjustable robust optimization problems \citep{bertsimas2022robust}. The actual solution, which we call here-and-now decision, shall be robust against some uncertain parameter under the premise that we can use wait-and-see decisions to react to the uncertainty once it reveals itself. For single-objective adjustable robust problems, there already exists a variety of solution methods. However, in practice, there is usually not only a single goal to be optimized but a variety of often contradicting objectives. In our example, we might not only want to minimize costs but also the carbon emissions of the energy supply. The arising multicriteria adjustable robust optimization (MARO) problems have not been studied in their general form in the literature yet. 

The concept of robust optimization is discussed in the book by \cite{bertsimas2022robust}. They pay special attention to the concept of adjustable robustness and give a good overview over the existing literature.
The concept of multicriteria optimization (MCO) is, for example, explained in the book by  \cite{ehrgott2005multicriteria}.
Multicriteria robust optimization (MRO) combines these two concepts. An overview over mathematical theory for MRO problems is given by \cite{ide2014concepts} and \cite{botte2019dominance}. One of the main challenges in computing solutions for MRO problems is that robustification and scalarization do not necessarily commute, as \cite{fliege2014robust} showed. 

MARO is a relatively new topic, that was first addressed mathematically by \cite{chuong2022adjustable}. They theoretically studied affinely adjustable problems, which can actually be interpreted as min-max robust (non-adjustable) problems, for which they considered point-based min-max robust efficiency. From the practical side, there is, for example, literature from energy system design \citep{yang2024two}, where a weighted sum approach is applied. \cite{groetzner2022multiobjective} consider a multicriteria regret approach, which is similar to an adjustable problem in that it can be formulated as a multi-level problem with three levels, which they can reduce to two levels. When it comes to MARO, one particularly complicating factor is decision uncertainty \citep{sinha2015solving}, meaning that at the time of the here-and-now-decision it is not clear which trade-off will be chosen for the choice of wait-and-see variables. These ideas were considered in the energy context by \cite{hollermann2021flexible}. Without decision uncertainty, problems can be simplified \citep{sinha2015towards}.

Stochastic optimization is another approach to deal with uncertainty. Multicriteria stochastic two stage problems are similar to MARO problems, but optimize the expected value instead of considering the maximum. In recent years, there were approaches in case studies to deal with such problems by applying the augmented $\varepsilon$-constraint method \citep{tiong2023two} or weighted sum inspired methods \citep{mena2023multi}.

In this paper, we propose and evaluate different definitions of efficient solutions to MARO problems under decision uncertainty. To find efficient solutions, we extend the ideas for MRO problems of \cite{ehrgott2014minmax} to multicriteria \emph{adjustable} robust problems and compare the computational concepts. For the comparison, we focus especially on whether efficiency concepts produce an interpretable image in objective space, that can easily be understood by decision makers, and on giving bounds that will hold for all scenarios. Moreover, we compare potential changes in the computational structure of the problem. We also examine relations between the different definitions.

Throughout the paper, we focus on problems for which we are able to solve the scalar adjustable robust problem for every objective function, as we explicitly exploit this structure in our computational approaches.

The remainder of this article is structured as follows. In Section~\ref{basics} we summarize the basics of MCO and MRO. In Section~\ref{Problem} we introduce the general MARO problem and provide definitions that characterize specific solutions for the problem. In Section~\ref{approaches}, we look at three different computational approaches, that can also be seen as further solution definitions based on solving scalar three-level problems. We examine their properties, evaluate whether they lead to meaningful solutions to the problem in the sense of Section~\ref{Problem} and compare them in Section~\ref{comparison}. In Section \ref{conclusionandoutlook} we summarize our findings and give an outlook on future research.

\section{Fundamentals} \label{basics}

In the following, we give a brief introduction to multicriteria and multicriteria robust optimization. The theory of these topics is naturally based on comparing values, vectors and sets.

\subsection{Multicriteria optimization}
In this section, we repeat some basics about MCO and introduce our notation. In MCO, we want to optimize different real valued objective functions $f_i$ at the same time. A parametric MCO problem with fixed parameter $u$ looks as follows. We want to solve for some set $\mathcal{X}(u)$ 
\begin{align}
	\tag{MUO}\label{MUO}
	\begin{split}
		\min_x &  \; f(x,u) \\
		\text{s.t.} & \; x \in \mathcal{X}(u),
	\end{split}
\end{align}
where we use for the sake of readability the notation
\begin{align*}
	f(x,u):= \begin{pmatrix}
		f_1(x, u) \\ ... \\f_n(x, u)
	\end{pmatrix}.
\end{align*}
Note that the parameter $u$ is not necessary for pure MCO problems and can therefore be omitted but will be used later to make the problem robust. If $f_1, ..., f_n$ do not share a common minimum, this problem needs a multicriteria understanding of optimality.

\begin{defi} \label{pref}
	For $y, y' \in \mathbb{R}^n$ we define three relations
	\begin{align*}
		y \leqq y' &:\Leftrightarrow y_i \leq y'_i &\forall i \in \{1,...,n\} \\
		y \leq y' &:\Leftrightarrow y \neq y', y_i \leq y'_i & \forall i \in \{1,...,n\} \\
		y < y' &:\Leftrightarrow y_i < y'_i & \forall i \in \{1,...,n\}.
	\end{align*}
    $\geqq, \geq$ and $>$ are defined analogously.
\end{defi}

With the notation
\begin{align*}
	\mathbb{R}^n_{[\geqq/\geq/>]} := \big\{ x \in \mathbb{R}^n \mid x \; [\geqq/\geq/>] \; 0_n \big\} ,
\end{align*}
where $0_n$ is the n-dimensional vector consisting only of zeros, we can equivalently write 
\begin{align*}
	y \;[\leqq / \leq / <] \; y' 
	\Leftrightarrow y \in y' - \mathbb{R}^n_{[\geqq/\geq/>]} 
  \Leftrightarrow y' \in y + \mathbb{R}^n_{[\geqq/\geq/>]}.
\end{align*}

Moreover, we define for two functions $f,g: A \rightarrow \mathbb{R}^n$, where $A$ is an arbitrary set,
\begin{align*}
	f \leqq g &:\Leftrightarrow \forall a \in A: f(a) \leqq g(a), \\
	f \leq g &:\Leftrightarrow (\forall a \in A: f(a) \leqq g(a))  \wedge (\exists a \in A: f(a) \leq g(a)). 
\end{align*}

Now, we can define the concepts of efficiency and nondominance as introduced in \citep[Definition 2.1]{ehrgott2005multicriteria}.

\begin{defi}
	A feasible solution $x \in \mathcal{X}$ is called \emph{[strictly/-/weakly] efficient} if
	\begin{align*}
		\nexists x' \in \mathcal{X}\backslash\{x\}: f(x') \; [\leqq/\leq/<] \; f(x).
	\end{align*} 
	If $x$ is [strictly/-/weakly] efficient, $f(x)$ is called \emph{[strictly/-/weakly] nondominated point}. If there is $x'$ such that $ f(x') \; [\leqq/\leq/<] \; f(x)$, we say $f(x')$ \emph{[weakly/-/strictly] dominates} $f(x)$. The set of all nondominated points is denoted $f(\mathcal{X})_N$ and the set of all efficient points is denoted $\mathcal{X}_E$. We sometimes use the same terminology for sets $S \subset \mathbb{R}^n$ without specifying a function. In this case, we refer to $id: S \rightarrow S$. Then, $S_{N^{min}}$ and $S_{N^{max}}$ correspond to the set of nondominated points with respect to $id$ and $-id$.
\end{defi}

There are different approaches to find [strictly/-/weakly] efficient solutions, like weighted sum scalarization, $\varepsilon$-constraint method and weighted Chebyshev scalarization. Scalarizations transform the multicriteria to a single criteria problem and guarantee that solutions to the single criteria problem are also efficient for the multicriteria one.

\subsection{Multicriteria Robust Optimization}
The formal robust counterpart of \eqref{MUO} is given by
\begin{align}
	\tag{MRO} \label{MRO}
	\begin{split}
		\min_x \max_{u}& \; f(x,u) \\
		\text{s.t.} & \; x \in \mathcal{X}(u) \; \forall u \in \mathcal{U}.
	\end{split}
\end{align}

It is not yet defined how to understand the inner maximization problem.
A lot of approaches to this problem can be found in \citep{botte2019dominance}, of which we will later refer to the following. For the sake of completeness, we give a strict and a weak version, whenever the definition is straightforward.
\begin{defi} \label{rmcodef}
	A feasible solution $x$ to \eqref{MRO} is 
	\begin{enumerate}
		\item \emph{{[strictly/-/weakly]} point-based robust min-max efficient} if
		\begin{align*} \nexists x' \in \mathcal{X} \backslash\{x\}: \begin{pmatrix}
				\max_{u} f_1(x', u) \\ ... \\ \max_{u} f_n(x', u)
			\end{pmatrix}
			\; [\leqq / \leq / <] \; \begin{pmatrix}
				\max_{u} f_1(x, u) \\ ... \\ \max_{u} f_n(x, u)
			\end{pmatrix}.
		\end{align*}
		\item \emph{{[strictly/-/weakly]} flimsy robust efficient} if 
		\begin{align*}
			\exists u \in \mathcal{U} \; \nexists x' \in \mathcal{X}\backslash\{x\}: f(x', u) \; [\leqq / \leq / <] \; f(x, u).
		\end{align*}
		\item \emph{{[strictly/-/weakly]} highly robust efficient} if
		\begin{align*}
			\forall u \in \mathcal{U} \; \nexists x' \in \mathcal{X}\backslash\{x\}: f(x',u) \;  [\leqq / \leq / <]\; f(x,u).
		\end{align*}
		\item \emph{{[strictly/-]} multi-scenario efficient} if 
		\begin{align*}
			\nexists x' \in \mathcal{X}\backslash\{x\}:f(x',\cdot) \;  [\leqq / \leq] \; f(x,\cdot).
		\end{align*}
	\end{enumerate}
\end{defi}

\subsection{Comparing Sets}

The following definitions will be needed in the upcoming Section~\ref{Problem}, where we will extend Definition~\ref{rmcodef} to set comparison. There are different ways to compare sets. The first two of the following are widely used in the literature and are referred to as [-/strict] upper and lower type set relation \citep{eichfelder2024set}.

\begin{defi} \label{preords}
	We define for $A, B \subset \mathbb{R}^n, \lambda \in \mathbb{R}^n_{\geq} $ the relations
	\begin{align*}
		&A \; [\leqq / <]^{u} \; B& &:\Leftrightarrow A \subset B - \mathbb{R}^n_{[\geqq/>]}  \\
		&A \; [\leqq / <]^{l} \; B& &: \Leftrightarrow B \subset A + \mathbb{R}^n_{[\geqq/ >]}  \\
		&A \; [\leqq / <]^{\lambda, min} \; B& &: \Leftrightarrow \min_{a \in A} \lambda^T a  \;[\leq / <] \; \min_{b \in B} \lambda^T b
	\end{align*}
	and write $A \; [\leqq / <]^* \; B$ to refer to any of them.
\end{defi}

\begin{rem}
	The $\leqq^*$-relations from Definition~\ref{preords} are preorders on the power set $\mathcal{P}(\mathbb{R}^n)$.
\end{rem}

\begin{rem} \label{cardone}
	If the cardinality of $A$ and $B$ is one, the relations $\leqq^*$ and  $<^*$ from Definition~\ref{preords} correspond for $* \in \{u,l\}$ to the relations $\leqq$ and $<$ from Definition~\ref{pref}. Moreover, we have for $A=\{a\}, B=\{b\},\; a,b \in \mathbb{R}^n, \lambda \in \mathbb{R}^n_{\geq} $ the implication
	\begin{align*}
		a \; [\leqq/<] \; b \; \Rightarrow \;  A \; [\leqq / <]^{\lambda, min} \; B.
	\end{align*}
\end{rem}

\section{Multicriteria Adjustable Robust Optimization Problems} \label{Problem}

In the following we consider the multicriteria adjustable robust optimization problem 
\begin{align}	\tag{MARO}\label{MARO}
	\begin{split}
		\min_x \max_{u} \min_y & \begin{pmatrix}
			f_1(x,u,y) \\ ... \\ f_n(x,u,y)
		\end{pmatrix} \\
		\text{s.t. }& x \in \mathcal{X} \subset \mathbb{R}^{n_x} \\
		& u \in \mathcal{U} \subset \mathbb{R}^{n_u} \\
		& y \in \mathcal{Y}(x, u) \subset \mathbb{R}^{n_y},
	\end{split}
\end{align}
where $x$ is the here-and-now decision, $u$ is the realization of uncertainty, and $y$ is the wait-and-see decision. We assume that $f_i$ is continuous for all $i \in \{1,...,n\}$ and $\mathcal{X,U}$ and $\mathcal{Y}(x,u)$ are compact sets. For the rest of the article, we will, for the sake of readability, sometimes drop the arguments of $\mathcal{Y}$ and write $f$ for the objective function vector. This problem combines the difficulties of MCO and robust/multi stage optimization and raises the question of what this term actually means.

\definecolor{lightgray}{RGB}{230, 230, 230}

\begin{figure}
	\centering
	\resizebox{0.6\textwidth}{!}{%
		\begin{tikzpicture}[scale=0.5, domain=2:10]
			\draw[very thin,color=lightgray] (-0.1,-0.1) grid (10.1,10.1);
			
			\draw[->,color=gray] (-0.1,0) -- (10.1,0) node[right] {$f_1$};
			\draw[->,color=gray] (0,-0.1) -- (0,10.1) node[above] {$f_2$};
			
			\fill[blue, opacity=0.2, domain=0.5:8.5] plot (\x, {10/((\x+1.5)) + 1.5}) -- (8,5) -- (4,8) -- cycle;
			\fill[blue, opacity=0.2, domain=0.5:8.5] plot (\x, {10/((\x+1.5)) + 2.5}) -- (8,6) -- (4,9) -- cycle;
			\fill[red, opacity=0.2, domain=0.5:6.5] plot (\x, {5/(0.25*(\x+1.5))-1.0}) -- (5,6) -- cycle;
			\fill[red, opacity=0.2, domain=0.5:6.5] plot (\x, {5/(0.25*(\x+1.5))}) -- (5,7) -- cycle;
			
			\draw[color=red, domain=0.5:6.5] plot (\x, {5/(0.25*(\x+1.5))-1.0}) node[below right] {$f(x_1, u_1, \mathcal{Y}(x_1, u_1))$};
			\draw[color=red, domain=0.5:6.5] plot (\x, {5/(0.25*(\x+1.5))}) node[below right] {$f(x_1, u_2, \mathcal{Y}(x_1, u_2))$};
			\draw[color=blue, domain=0.5:8.5] plot (\x, {10/(\x+1.5) + 1.5}) node[above] {$f(x_2, u_1, \mathcal{Y}(x_2, u_1))$};
			\draw[color=blue, domain=0.5:8.5] plot (\x, {10/(\x+1.5) + 2.5}) node[above] {$f(x_2, u_2, \mathcal{Y}(x_2, u_2))$};
		\end{tikzpicture}
	}
	\caption{Wait-and-see options for two different first stage decisions $x_1,x_2$ and two different uncertainty realizations $u_1,u_2$. The areas refer to the named sets, and the respective sets of nondominated points are highlighted as lines.}\label{visualization} 
\end{figure}
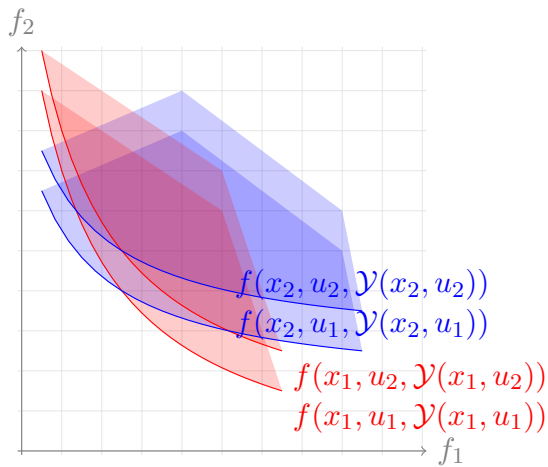

Figure~\ref{visualization} shows the setting for discrete sets $\mathcal{X}=\{x_1, x_2\}, \;  \mathcal{U}=\{u_1, u_2\}$ and continuous $\mathcal{Y}$. Under decision uncertainty - that is, if the choice of $y$ is not predetermined by a given a priori preference - it is natural to reduce the inner stage of the problem to the efficient set $\mathcal{Y}(x,u)_{E}$, that minimizes $f(x,u,y)$ for fixed $x$ and $u$. There exists now a variety of possibilities of how to understand this problem. In the following definition, we extend MRO to include set comparison at the inner level.

\begin{defi} \label{marodef}
	We call $x \in \mathcal{X}$ for $[\leqq/<]^*$
	\begin{enumerate}
		\item \emph{{[strictly/weakly]} flimsy MARO-efficient} if 
		\begin{align*}
			\exists u \in \mathcal{U}  \;  \nexists x' \in \mathcal{X}\backslash\{x\} : f(x',u, \mathcal{Y}(x',u)_{E}) \;  [\leqq / <]^*  \; f(x, u, \mathcal{Y}(x,u)_{E}),
		\end{align*}
		\item \emph{{[strictly/weakly]} highly MARO-efficient} if 
		\begin{align*}
			\forall u \in \mathcal{U}  \;  \nexists x' \in \mathcal{X}\backslash\{x\}: f(x',u, \mathcal{Y}(x',u)_{E}) \;  [\leqq / <]^*  \; f(x, u, \mathcal{Y}(x,u)_{E}),
		\end{align*}
		\item \emph{strictly multi-scenario MARO-efficient} if  
		\begin{align*}
			\nexists x' \in \mathcal{X}\backslash\{x\}: f(x',\cdot, \mathcal{Y}(x',\cdot)_{E}) \;  \leqq ^* \;  f(x,\cdot, \mathcal{Y}(x,\cdot)_{E}),
		\end{align*}
		where the relation is a generalization of the comparison of vector valued functions to set-valued functions, that is defined for $f,g: A \rightarrow \mathcal{P}(\mathbb{R}^n)$, where $A$ is an arbitrary set,  as $f \leqq^* g :\Leftrightarrow \forall a \in A: f(a) \leqq^* g(a)$.
	\end{enumerate}
\end{defi}

There are several easy to observe implications between the concepts.

\begin{rem}
	Every feasible solution $x$ to \eqref{MARO} that is strictly {[flimsy/ highly]} MARO-efficient is weakly {[flimsy/highly]} MARO-efficient.
	Moreover, we observe that every {[strictly/weakly]} highly MARO-efficienct $x$ is [strictly/weakly] flimsy MARO-efficienct, and every strictly highly MARO-efficienct $x$ is strictly multi-scenario MARO-efficienct. 
\end{rem}

Weakly flimsy MARO-efficient solutions for $<^l$ can be found by considering the MCO
\begin{align*}
	\min_{x \in \mathcal{X},u \in \mathcal{U},y \in \mathcal{Y}} f(x,u,y)
\end{align*}
and are closer related to best cases than to robustness against worst cases.
{[Strictly/Weakly]} highly MARO-efficient solutions are not very likely to exist in practice, as their existence would mean that there is a common optimal $x$ for all elements of the uncertainty set.

We have the following coherence results.

\begin{rem}
	If $\mathcal{U}$ is a singleton, i.e. $\mathcal{U} = \{u\}$, the problem can be reduced to a single stage MCO problem, and we would require in MCO without decision uncertainty
	\begin{align*}
		\exists y \in \mathcal{Y}(x) \; \nexists x' \in \mathcal{X}\backslash\{x\}, y' \in \mathcal{Y}(x'): f(x',u,y') \leq f(x,u,y),
	\end{align*}
	for $x \in \mathcal{X}$ to be part of an efficient solution. However, for all concepts from Definition~\ref{marodef} we generalize the efficiency understanding in order to take the wait-and-see-variables into account and only require 
	\begin{align*}
		\nexists x' \in \mathcal{X}\backslash\{x\}: f(x,u, \mathcal{Y}(x,u)) \; [\leqq/<]^* \; f(x',u,\mathcal{Y}(x',u)).
	\end{align*}
	This corresponds to an understanding of the inner problem as a set optimization problem. Therefore, our definitions are coherent with MCO under decision uncertainty.
\end{rem}

\begin{lem}
	If $\mathcal{Y}$ is a singleton, the MARO problem reduces to an MRO problem. For $x \in \mathcal{X}$ and $[\leqq/<]^*, *\in \{u,l\}$ we have
	\begin{align*}
		&x  \text{ [strictly/weakly] flimsy/highly robust efficient}\\
		&\Leftrightarrow x \text{ [strictly/weakly] flimsy/highly MARO-efficient} 
	\end{align*}
	and 
	\begin{align*}
		& x  \text{ strictly multi-scenario efficient}\\
		&\Leftrightarrow x \text{ strictly multi-scenario MARO-efficient} .
	\end{align*}	
	In case of $[\leqq / <]^{\lambda, min}$ for some $\lambda \in \mathbb{R}^n_{\geq}$, only the implications hold but not the reverse directions.
\end{lem}

\begin{pf}
	This follows by applying Remark~\ref{cardone} to the definitions.
\end{pf}

\begin{figure}
	\parbox{0.5\textwidth}{
		\begin{tikzpicture}[scale=0.5, domain=0:10]
			\draw[very thin,color=lightgray] (-0.1,-0.1) grid (10.1,10.1);
			
			\draw[->,color=gray] (-0.1,0) -- (10.1,0) node[right] {$f_1$};
			\draw[->,color=gray] (0,-0.1) -- (0,10.1) node[above] {$f_2$};
			
			\fill[color=black] (3, 7) circle (3pt) node[above] {$f(x_1, u_1,y)$};
			\fill[color=black] (5, 5) circle (3pt) node[above right] {$f(x_1, u_2,y)$};
			\fill[color=black] (8, 4) circle (3pt) node[above] {$f(x_1, u_3,y)$};
			\fill[color=blue] (2, 6) circle (3pt) node[below] {$f(x_2, u_3,y)$};
			\fill[color=blue] (4, 4) circle (3pt) node[below] {$f(x_2, u_2,y)$};
			\fill[color=blue] (7, 3) circle (3pt) node[below] {$f(x_2, u_1,y)$};
		\end{tikzpicture}
	}
	\begin{minipage}{0.5\textwidth}
		\begin{tikzpicture}[scale=0.5, domain=0:10] 
			\draw[very thin,color=lightgray] (-0.1,-0.1) grid (10.1,10.1);
			
			\draw[->,color=gray] (-0.1,0) -- (10.1,0) node[right] {$f_1$};
			\draw[->,color=gray] (0,-0.1) -- (0,10.1) node[above] {$f_2$};
			
			\fill[color=black] (3, 7) circle (3pt) node[above] {$f(x_1, u_1,y)$};
			\fill[color=black] (5, 5) circle (3pt) node[above] {$f(x_1, u_2,y)$};
			\fill[color=blue] (2, 2) circle (3pt) node[below] {$f(x_2, u_1,y)$};
			\fill[color=blue] (4, 4) circle (3pt) node[below] {$f(x_2, u_2,y)$};
		\end{tikzpicture}
	\end{minipage}
	\caption{In the left picture, $x_1$ is strictly multi-scenario MARO-efficient for $\leqq^l$ but not efficient for \eqref{SMARO}, and in the right picture the situation is vice versa. Moreover, it is obvious that in the right picture $x_1$ should not be called efficient for any reasonable efficiency concept. Also, in the left picture, $x_1$ seems worse than $x_2$ at the first glance. However, in scenarios $u_1$ and $u_3$ it is impossible to say which $x$ is best. Therefore, it seems reasonable to call $x_1$ and $x_2$ both strictly multi-scenario MARO-efficient for $\leqq^l$. In both pictures, we assume $\mathcal{Y}=\{y\}$ and therefore we are actually in a non-adjustable robust setting.}
	\label{Paretofrontandcounterexample}
\end{figure}
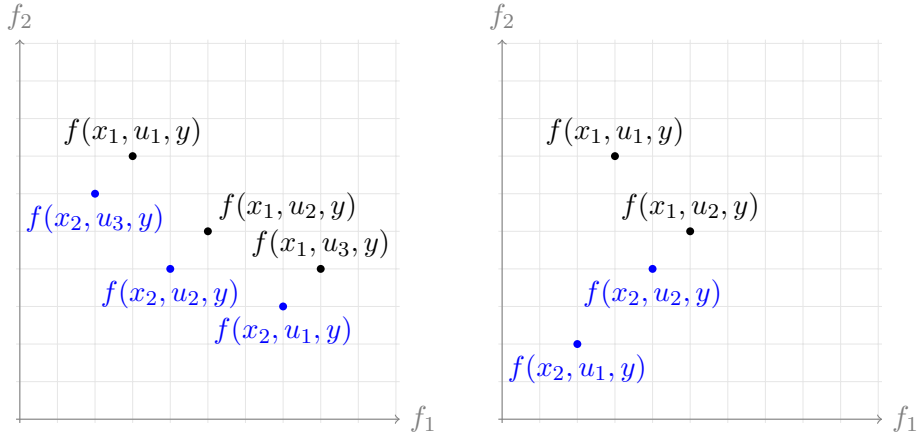

We suggested the above definitions of efficiency for \eqref{MARO} because they align well with our intuition of what a solution should look like.
In contrast, the following definition seems natural but unfortunately causes counterintuitive results: Assume we consider the problem as a set optimization problem at every stage, where we call $x$ that corresponds to points in

\begin{align*}
	\tag{SMARO}\label{SMARO}
		&\left\{ \bigcup_{x \in \mathcal{X}} \left\{ \bigcup_{u \in \mathcal{U}} \left\{ \bigcup_{y \in \mathcal{Y}(x,u)}   \begin{pmatrix} f_1(x,u,y) \\ ... \\ f_n(x,u,y) \end{pmatrix} \right\}_{N^{min}} \right\}_{N^{max}} \right\}_{N^{min}}
\end{align*}
efficient for \eqref{SMARO}.
Then, Figure~\ref{Paretofrontandcounterexample} shows a situation in which this understanding of efficiency leads to undesirable efficient solutions.

\section{Computational Approaches to MARO} \label{approaches}
For the definitions of MARO-efficiency considered so far, it is not immediately clear how to find solutions that fulfill them. We will now study three different approaches, that can be interpreted as further MARO-efficiency concepts, in more detail. These efficiency concepts share the advantage that it follows immediately from the definition how to find efficient solutions for them. 
In the exploration of these approaches, we focus on the following properties, as they are relevant for determining the most suitable approach:
\begin{itemize}
	\item Are there bounds for efficient solutions over all scenarios?
	\item Do the concepts lead to an interpretable image in objective space? Does it resemble a common multicriteria optimization Pareto front in the sense that it only contains nondominated or at least weakly nondominated points? This property is especially valuable for decision makers, as the decision-making process can then be carried out similarly as for standard (non-adjustable robust) MCO problems. 
	\item Do we preserve the problem's inherent structure when scalarizing the problem? For example, if we solve a multicriteria flow problem at the inner level, will our concept also solve a flow problem?
	\item Are solutions that are efficient according to the respective definition also efficient for some definition from Section~\ref{Problem}? 
\end{itemize}

All considered concepts rely on the solution of a scalar three stage problem. There is no standard procedure to solve these problems for general problem structures. An algorithm for such problems is given by \cite{zeng2013solving} for discrete or polynomial uncertainty set and continuous, linear inner stage and by \cite{yue2019projection} for discrete or polynomial uncertainty set and mixed integer linear inner stage. In both cases, the general structure is a linear problem with right hand side uncertainty. These algorithms use duality to get rid of the inner stage and then apply adaptive discretization \citep{blankenship1976infinitely}. This requires computing the global optimal solution of the middle stage problem.

\subsection{Weighted sum efficiency}
In this section, we look at an approach that scalarizes first. The idea is inspired by the weighted sum scalarization approach for MCO \citep[Proposition 3.9]{ehrgott2005multicriteria}.

We consider the function $f_{\lambda}: \mathcal{X} \rightarrow \mathbb{R}$,
\begin{align*}
	f_{\lambda}(x) := &\max_{u \in \mathcal{U}} \min_{y \in \mathcal{Y}} \sum_{i=1}^n \lambda_i f_i(x,u,y),
\end{align*}
where $\lambda \in \Lambda := \{ \lambda \in \mathbb{R}^n_{\geq} \mid \sum_{i=1}^n \lambda_i = 1\}$
and introduce the notation
\begin{align*}
	\mathcal{U}_{\lambda}^*(x) = \argmax_{u \in \mathcal{U}} \min_{y \in \mathcal{Y}} \sum_{i=1}^n \lambda_i f_i(x,u,y), \\ \mathcal{Y}_{\lambda}^*(x,u) = \argmin_{y \in \mathcal{Y}} \sum_{i=1}^n \lambda_i f_i(x,u,y).
\end{align*}

\begin{defi}
	We call $x \in \mathcal{X}$ \emph{[strictly/-] weighted sum MARO-efficient} for $\lambda \in \Lambda$ with guarantee $\mathcal{G}_\lambda := f_\lambda(x)$ if
	\begin{align*}
		\nexists x' \in \mathcal{X} \backslash\{x\}: f_{\lambda}(x') \; [\leq /<] \; f_{\lambda}(x).
	\end{align*}
\end{defi}

Weighted sum MARO-efficient solutions can be found by solving a three stage optimization problem. Moreover, when implementing a weighted sum MARO-efficient solution, we immediately get a bound which holds independent of the scenario.

\begin{rem} \label{weightedbound}
	For a weighted sum MARO-efficient solution $x \in \mathcal{X}$ for $\lambda \in \Lambda$ we obtain the bound
	\begin{align*}
		\forall u \; \exists y: \; \sum_{i=1}^n \lambda_i f_i(x,u,y) \leq \mathcal{G}_\lambda.
	\end{align*}
\end{rem}

We now look at the image of weighted sum MARO-efficient solutions in the objective space, as decision-makers usually decide in this space.
 
\begin{defi}
	For $\lambda \in \Lambda$ we define 
	\begin{align*}
		im^f(\lambda) := \{ f(x, u, y) \mid &x \text{ weighted sum MARO-efficient for } \lambda, \\ & u \in \mathcal{U}_{\lambda}^*(x), y \in \mathcal{Y}_{\lambda}^*(x,u)\}
	\end{align*} 
 	and the weighted sum MARO-efficiency objective space image
	\begin{align*}
		im^f(\Lambda) := \bigcup_{\lambda \in \Lambda}im^f(\lambda).
	\end{align*}
\end{defi}

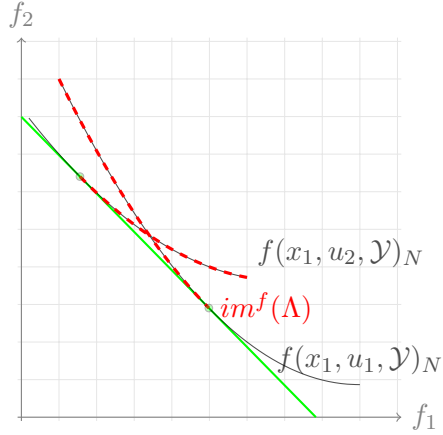
\begin{figure}
	\centering
		\begin{tikzpicture}[scale=0.5, domain=0:10] 
			\draw[very thin,color=lightgray] (-0.1,-0.1) grid (10.1,10.1);
			
			\draw[->,color=gray] (-0.1,0) -- (10.1,0) node[right] {$f_1$};
			\draw[->,color=gray] (0,-0.1) -- (0,10.1) node[above] {$f_2$};
			\draw[domain=0:7.83078,smooth,variable=\x, color=green, thick=1] plot ({\x},{-1.02051*\x + 7.99139}) node[above left] {};
			\draw[domain=1:9,smooth,variable=\x, opacity=0.7, thin] plot ({\x},{23/180*\x*\x-413/180*\x+67/6}) node[above] {$f(x_1, u_1, \mathcal{Y})_N$};
			\draw[domain=0.2:6,smooth,variable=\x, opacity=0.7, thin] plot ({\x},{3/32*\x*\x-21/16*\x+263/32}) node[above right] {$f(x_1, u_2, \mathcal{Y})_N$};
			\draw[domain=1:4.98498,smooth,variable=\x, color=red, very thick, dashed] plot ({\x},{23/180*\x*\x-413/180*\x+67/6}) node[right] {$im^f(\Lambda)$};
			\draw[domain=1.5573:6,smooth,variable=\x, color=red, very thick, dashed] plot ({\x},{3/32*\x*\x-21/16*\x+263/32}) node[above right] {};
			\draw[fill=green, opacity=0.3] (4.98498,-1.02051*4.98498 + 7.99139) circle (0.1);
			\draw[fill=green, opacity=0.3] (1.5573,-1.02051*1.5573 + 7.99139) circle (0.1);
		\end{tikzpicture}
	\caption{$im^f(\Lambda)$ is shown in dashed red. Not all points in $im^f(\Lambda)$ are weakly nondominated, and for the $\lambda$ that is associated with the green line, $im^f(\lambda)$ is a disconnected set.}\label{pintersect}
\end{figure}

\begin{rem} \label{lambdafront}
	There are several drawbacks of weighted sum MARO-efficiency. Figure~\ref{pintersect} depicts some undesired features of such a solution. Not all points in $im^f(\Lambda)$ are weakly nondominated. Further, it is open what to do if $im^f(\lambda)$ is a disconnected set for $\lambda \in \Lambda$, as for the $\lambda$ that is associated with the green line.
\end{rem}

Finally, we can relate weighted sum MARO-efficiency to MARO-efficiency definitions from Section~\ref{Problem}.

\begin{thm} \label{firstweightedlem}
	Every for $\lambda \in \Lambda$ strictly weighted sum MARO-efficient solution is strictly multi-scenario MARO-efficient with respect to the preorder $\leqq^{\lambda, min}$.
\end{thm}

\begin{pf}
	Let $x$ be strictly weighted sum MARO-efficient. This is defined as
	\begin{align*}
		\nexists x' \in \mathcal{X}\backslash \{x\}: \max_u \min_{y \in \mathcal{Y}} \sum_i \lambda_i f_i(x',u,y) \leq \max_u \min_{y \in \mathcal{Y}} \sum_i \lambda_i f_i(x,u,y).
	\end{align*}
	By contraposition we obtain
	\begin{align*}
		\nexists x' \in \mathcal{X}\backslash \{x\} \; \forall u \in \mathcal{U}: \min_{y \in 	\mathcal{Y}} \sum_{i} \lambda_i f_i(x', u, y) \leq \min_{y \in \mathcal{Y}} \sum_i \lambda_i f_i(x,u,y).
	\end{align*}
	Note that, exploiting \citep[Theorem 2.21]{ehrgott2005multicriteria}, $\min_{y \in \mathcal{Y}}$ can be replaced by $\min_{y \in \mathcal{Y}_E}$. Inserting the definition of $\leqq^{\lambda,min}$, we get
	\begin{align*}
		\nexists x' \in \mathcal{X}\backslash \{x\} \; \forall u \in \mathcal{U}: f(x', u, \mathcal{Y}(x',u)_{E}) \leqq^{\lambda,min} f(x,u, \mathcal{Y}(x,u)_{E}),
	\end{align*}
	which is by definition of the comparison of set-valued functions equivalent to
	\begin{align*}
		\nexists x' \in \mathcal{X}\backslash \{x\}: f(x',\cdot,\mathcal{Y}(x',\cdot)_{E}) \leqq^{\lambda,min} f(x,\cdot, \mathcal{Y}(x,\cdot)_{E}).
	\end{align*}
	This is precisely the definition of $x$ being strictly multi-scenario MARO-efficient for the preorder $\leqq^{\lambda, min}$. \qed
\end{pf}

\subsection{Constraint efficiency}

In this section, we look at an approach that also scalarizes first but is inspired by the $\varepsilon$-constraint method for MCO \citep[Chapter 4.1]{ehrgott2005multicriteria}. We set hard limits on all but one objective function and can therefore reduce the problem to a single criteria problem.
We consider for generating bound $\varepsilon \in \mathbb{R}^n$ and $j \in \{1,...,n\}$ the function $f_{\varepsilon,j}: \mathcal{X} \rightarrow \mathbb{R}$,
\begin{align*}
	f_{\varepsilon,j}(x) := &\max_{u \in \mathcal{U}} \left\{\min_{y \in \mathcal{Y}} f_j(x,u,y) \;  s.t. \; f_i(x,u,y) \leq \varepsilon_i \; \forall i \neq j\right\}.
\end{align*}

We follow the convention that minimization over an empty set has result $+\infty$ and maximization over an empty set has result $-\infty$.

\begin{defi}
	We call $x \in \mathcal{X}$ \emph{[strictly/-] constraint MARO-efficient} for the generating bound $\varepsilon \in \mathbb{R}^n$ and $j \in \{1,...,n\}$ with guarantee $\mathcal{G}_{\varepsilon,j}:= f_{\varepsilon,j}(x)$ if there is no $x'$ with 
	\begin{align*}
		f_{\varepsilon,j}(x') \; [</\leq] \; f_{\varepsilon,j}(x).
	\end{align*}
\end{defi}

The computation of constraint MARO-efficient solutions is again straight forward by solving a three stage problem.

\begin{thm} \label{switchlem}
	Let $x$ be strictly constraint MARO-efficient for $\varepsilon \in \mathbb{R}^n$ and $j \in \{1,...,n\}$. Then, $x$ is strictly constraint MARO-efficient for $\varepsilon' = (\varepsilon_1,..., \varepsilon_{j-1}, f_{\varepsilon,j}(x), \varepsilon_{j+1}, ..., \varepsilon_n)^T$, independent of the choice of $i \in \{1,...,n\}$ as the index of the chosen objective.
\end{thm}

\begin{pf}
	Let $x$ be strictly constraint MARO-efficient for $\varepsilon$ and $j$. Let $i \in\{1,...,n\}, i\neq j$. Assume there is $x' \in \mathcal{X} \backslash\{x\}$ such that $f_{\varepsilon',i}(x') \leq f_{\varepsilon',i}(x)$. 
	From 
 \begin{align*}
     		f_{\varepsilon,j}(x) = \max_{u \in \mathcal{U}} \left\{\min_{y \in \mathcal{Y}} f_j(x,u,y) \;  s.t. \; f_i(x,u,y) \leq \varepsilon_i \; \forall i \neq j\right\}
 \end{align*}we know that $\forall u \in \mathcal{U} \; \exists y \in \mathcal{Y}: f(x,u,y) \leq \varepsilon'$, i.e. $f_{\varepsilon',i}(x) \leq \varepsilon'_i$. Hence, due to our assumption, $\forall u \in \mathcal{U} \; \exists y \in \mathcal{Y}(x',u): f(x',u,y) \leq \varepsilon'$.
	By definition of $\varepsilon'$, this is a contradiction to $x$ being strictly constraint MARO-efficient for $\varepsilon$ and $j$.
	\qed
\end{pf}

Therefore, strict constraint MARO-efficiency is a well-defined approach, and the property is actually independent of $j$. However, in the following, we will see that also non-strictly constraint MARO-efficient solutions share some valuable properties. For example, they immediately come with good bounds for all scenarios and objectives.

\begin{rem} \label{constraintbound}
	For a [strictly/-] constraint MARO-efficient solution $x$ with respect to $j \in \{1,...,n\}$ and generating bound $\varepsilon \in \mathbb{R}^n$ we have for all $u \in \mathcal{U}$
	\begin{align*}
		\exists y \in \mathcal{Y}(x,u): (f_i(x,u,y) \leq \varepsilon_i \; \forall i \in \{1,...,n\}, i \neq j) \wedge (f_j(x,u,y) \leq \mathcal{G}_{\varepsilon,j}).
	\end{align*} 
\end{rem}

For decision-makers it can be very valuable to look at the results in objective space.

\begin{defi}
	We define for $j \in \{1,...,n\}, \epsilon \in \mathbb{R}^n$
	\begin{align*}
		im_j^f(\varepsilon):= \left(
			\varepsilon_1, ..., \varepsilon_{j-1}, \min_x f_{\varepsilon,j}(x), \varepsilon_{j+1}, ..., \varepsilon_n
		\right)^T
	\end{align*}
	and the constraint MARO-efficiency objective space image for a chosen constraint set $E \subset \mathbb{R}^n$
	\begin{align*}
		im_j^f(E):= \bigcup_{\varepsilon \in E} im_j^f(\varepsilon).
	\end{align*}
\end{defi}

\begin{lem} \label{constraintfront}
	For $E \subset \mathbb{R}^n$, all points in $im_j^f(E)$ are weakly nondominated.
\end{lem}

\begin{pf}
	Let $\varepsilon_1 \in E$ arbitrary but fixed. We will show that $im_j^f(\varepsilon_1)$ is weakly nondominated. Let $\varepsilon_2 \neq \varepsilon_1$ be any other point in $E$. If not $im_j^f(\varepsilon_2) < im_j^f(\varepsilon_1)$ in $\mathbb{R}^{n-1}$ where the $j$th entry of the vector is dropped, we are done. So let $\varepsilon_2 < \varepsilon_1$ without the $j$th vector entry. We then observe
	\begin{align*}
		im_j^f(\varepsilon_1)_j = \min_x f_{\varepsilon_1,j}(x) \leq \min_x f_{\varepsilon_2,j}(x) = im_j^f(\varepsilon_2)_j
	\end{align*}
	by definition of the optimization problems. Hence, $\varepsilon_1$ is weakly nondominated.
	\qed
\end{pf}
This result implies that $im_j^f(E)$ has a shape that is easy to understand for decision makers.

Finally, we put strict constraint MARO-efficiency in the context of Section~\ref{Problem}.

\begin{thm} \label{theostrictconst}
	If $x$ is strictly constraint MARO-efficient for $\varepsilon \in \mathbb{R}^n, j \in \{1,...,n\}$, it is strictly multi-scenario MARO-efficient for the preorder $\leqq^{l}$.
\end{thm}

\begin{pf}	We start by reformulating the definition. Recall that a point $x$ is strictly constraint MARO-efficient for $\varepsilon \in \mathbb{R}^n, j \in \{1,...,n\}$ if 
\begin{align*}
	\nexists x' \in \mathcal{X} \backslash\{x\}: f_{\varepsilon,j}(x') \leq f_{\varepsilon,j}(x)
\end{align*}
or equivalently
\begin{align*}
	\forall x' \in \mathcal{X} \backslash \{x\}: f_{\varepsilon,j}(x') > f_{\varepsilon,j}(x).
\end{align*}
	Let $x \in \mathcal{X}$ have this property. Then, by inserting the definition of $f_{\varepsilon,j}$ and taking a maximal $u$ from the left side, we obtain 
\begin{align*}
	\forall x' \in \mathcal{X}\backslash \{x\} \; \exists u \in \mathcal{U}: 
	&\min_{y \in \mathcal{Y}} f_j(x',u,y) \text{ s.t. $f_i(x',u,y) \leq \varepsilon_i \; \forall i \neq j$ } >  \\ 
	&\min_{y \in \mathcal{Y}} f_j(x,u,y) \text{ s.t. $f_i(x,u,y) \leq \varepsilon_i \; \forall i \neq j$ }.
\end{align*}
	Exploiting \citep[Theorem 2.21]{ehrgott2005multicriteria}, we can replace $\min_{y \in \mathcal{Y}}$ by $\min_{y \in \mathcal{Y}_E}$.
	Taking a minimal $y$ from the right side, we obtain with the definition 
	$
		\mathcal{Y}|(x,u,j):= \left\{y \in \mathcal{Y}(x,u)_{E}  \mid f_i(x,u, y) \leq \varepsilon_i \; \forall i \neq j\right\}
	$
\begin{align*}
		\forall x' \in \mathcal{X}\backslash \{x\} \; \exists u \in \mathcal{U} \; &\exists y \in \mathcal{Y}|(x,u,j)\; \forall y' \in \mathcal{Y}|(x',u,j):  \\
		&f_j(x',u,y') > f_j(x,u,y).
\end{align*}
	By definition of the relation $\leq$ and using the notation $\nleq$ to indicate that two elements are not in relation with respect to it, we can conclude
\begin{align*}
		\forall x' \in \mathcal{X}\backslash \{x\} \; \exists u \in \mathcal{U} \; &\exists y \in \mathcal{Y}|(x,u,j)\; \forall y' \in \mathcal{Y}|(x',u,j): \\& f(x',u,y') \nleq f(x,u,y).
\end{align*}
	For all $y \in \mathcal{Y}(x',u)_{E}$ for which $f_i(x',u,y') \leq \varepsilon_i \; \forall i \neq j$ does not hold, the statement is already clear and therefore 
\begin{align*}
		\forall x' \in \mathcal{X}\backslash \{x\} \; \exists u \in \mathcal{U} \; & \exists y \in \mathcal{Y}|(x,u,j)\; \\ & \forall y' \in \mathcal{Y}(x',u)_{E}: f(x',u,y') \nleq f(x,u,y).
\end{align*}
	By again negating the statement, we get
\begin{align*}
		\nexists x' \in \mathcal{X}\backslash\{x\} \; \forall u \in \mathcal{U} \;&\forall y \in \mathcal{Y}|(x,u,j)\; \\ & \exists y' \in \mathcal{Y}(x',u)_{E} : f(x',u,y') \leq f(x, u,y).
\end{align*}
	If the statement does not hold for all elements of a subset, it will also not hold for all elements in the whole set. Hence, 
\begin{align*}
	\nexists x' \in \mathcal{X}\backslash\{x\} \; \forall u \in \mathcal{U} \;&\forall y \in \mathcal{Y}(x,u)_{E} \; \exists y' \in \mathcal{Y}(x',u)_{E}: \\ &f(x',u,y') \leq f(x, u,y)
\end{align*}
	and this is by definition of $\leqq^l$ equivalent to
\begin{align*}
	\nexists x' \in \mathcal{X} \backslash\{x\}:f(x',\cdot, \mathcal{Y}(x',\cdot)_{E}) \leqq^{l} f(x,\cdot, \mathcal{Y}(x,\cdot)_{E}).
\end{align*}
	This is precisely the definition of $x$ being strictly multi-scenario MARO-efficient for the preorder $\leqq^{l}$. \qed 
\end{pf}

\subsection{Point-based min-max-min MARO-efficiency}
As an extension to point-based min-max robust efficiency, the following might be the most straight-forward approach.
We now look at
\begin{align*} 
	f^{pb}: \mathcal{X} \rightarrow \mathbb{R}^n, f^{pb}(x) := &\begin{pmatrix}
		\max_{u \in \mathcal{U}} \min_{y \in \mathcal{Y}} f_1(x,u,y) \\ ... \\ \max_{u \in \mathcal{U}} \min_{y \in \mathcal{Y}} f_n(x,u,y)
	\end{pmatrix} .
\end{align*}

\begin{defi}
	We call $x \in \mathcal{X}$ \emph{[strictly/-/weakly] point-based min-max-min MARO-efficient} if there is no $x' \in \mathcal{X} \backslash \{x\}: f^{pb}(x') \; [\leqq / \leq / <] \; f^{pb}(x)$.
\end{defi}

Weighted sum MARO-efficiency and constraint MARO-efficiency are based on a first-scalarize-then-robustify approach. Finding point-based min-max-min MARO-efficient solutions follows a first-robustify approach. If the resulting MCO problem is solved with a scalarization technique, a first-robustify-then-scalarize approach is used.

In the following, we observe that for point-based min-max-min MARO-efficient solutions we do not get any bounds, except for trivial ones.

\begin{defi}
	For a compact set $A \subset \mathbb{R}^n$, we define its \emph{ideal point} with respect to minimization and maximization as
	\begin{align*}
		A_{I^{min}}:= \begin{pmatrix}
			\min_{a \in A} a_1 \\... \\ \min_{a \in A} a_n
		\end{pmatrix},A_{I^{max}}:= \begin{pmatrix}
			\max_{a \in A} a_1 \\... \\ \max_{a \in A} a_n
		\end{pmatrix}.
	\end{align*}
\end{defi}

\begin{figure}
	\centering
	\begin{tikzpicture}[scale=0.5, domain=0:10]
		\draw[very thin,color=lightgray] (-0.1,-0.1) grid (10.1,10.1);
		
		\draw[->,color=gray] (-0.1,0) -- (10.1,0) node[right] {$f_1$};
		\draw[->,color=gray] (0,-0.1) -- (0,10.1) node[above] {$f_2$};
		
		\draw plot[smooth] coordinates {(1, 9) (4,6)} node[above right] {$f(x_1, u_1, \mathcal{Y})=f(x_1, u_2, \mathcal{Y})$};
		\fill[color=black] (5,4) circle (3pt) node[below] {$f(x_2, u_1, \mathcal{Y})$};
		\fill[color=black] (8,2) circle (3pt) node[below] {$f(x_2, u_2, \mathcal{Y})$};
		\fill[color=blue] (1,6) circle (3pt) node[below] {$f^{pb}(x_1)$};
		\fill[color=blue] (8,4) circle (3pt) node[above] {$f^{pb}(x_2)$};
	\end{tikzpicture}
	\caption{ $f^{pb}(x_1)$ corresponds to $f(x_1, \mathcal{U},\mathcal{Y})_{I^{min}}$ and $f^{pb}(x_2)$ corresponds to $f(x_2, \mathcal{U},\mathcal{Y})_{I^{max}}$, and $x_1$ as well as $x_2$ are point-based robust min-max-min MARO-efficient.}\label{badexample}
\end{figure}
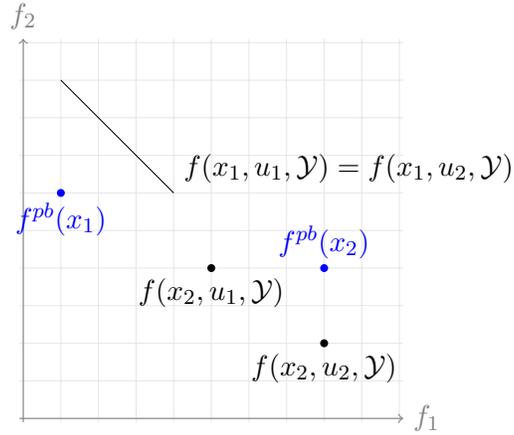

\begin{rem} \label{pointest}
	If $x$ is point-based min-max-min MARO-efficient, we have with the definition $f(x, \mathcal{U}, \mathcal{Y}):= \bigcup_{u \in \mathcal{U}}f(x,u, \mathcal{Y}(x,u))$
	\begin{align*}
		f(x,\mathcal{U},\mathcal{Y})_{I^{min}} \leq f^{pb}(x) \leq f(x,\mathcal{U},\mathcal{Y})_{I^{max}}
	\end{align*}
	and in general there is no better approximation as Figure~\ref{badexample} shows. 
\end{rem}

\begin{figure}
	\centering
	\begin{tikzpicture}[scale=0.5, domain=2:10]
		\draw[very thin,color=lightgray] (-0.1,-0.1) grid (10.1,10.1);
		
		\draw[->,color=gray] (-0.1,0) -- (10.1,0) node[right] {$f_1$};
		\draw[->,color=gray] (0,-0.1) -- (0,10.1) node[above] {$f_2$};
		
		\draw plot[smooth] coordinates {(1+2, 5+1) (1.5+2,3.5+1) (3+2,3+1)} node[above right] {$f(x_2, u_1, \mathcal{Y})_N$};
		\draw plot[smooth] coordinates {(2+2, 4+1) (2.5+2,2.5+1) (4+2,2+1)} node[above right] {$f(x_2, u_2, \mathcal{Y})_N$};
		\draw plot[smooth] coordinates {(7,2) (8,1)} node[right] {$f(x_3, u_1, \mathcal{Y})_N$};
		\draw plot[smooth] coordinates {(8,3) (9,2)} node[right] {$f(x_3, u_2, \mathcal{Y})_N$};
		\draw plot[smooth] coordinates {(0, 10) (1,9)} node[above right] {$f(x_1, u_1, \mathcal{Y})_N$};
		\draw plot[smooth] coordinates {(2,8) (3,7)} node[right] {$f(x_1, u_2, \mathcal{Y})_N$};
		\fill[color=blue] (2+2,3+1) circle (3pt) node[below left] {$f^{pb}(x_2)$};
		\fill[color=blue] (8,2) circle (3pt) node[above left] {$f^{pb}(x_3)$};
		\fill[color=blue] (2,9) circle (3pt) node[below right] {$f^{pb}(x_1)$};
	\end{tikzpicture}
	\caption{Set of nondominated points for the inner stage for three different first stage decisions $x_1, x_2, x_3$ and two different uncertainty realizations $u_1,u_2$. $f^{pb}(x_1)$ is weakly dominated by points in $f(x_1,\mathcal{U}, \mathcal{Y})$, while there are no points in $f(x_2,\mathcal{U}, \mathcal{Y})$ that weakly dominate $f^{pb}(x_2)$ and $f^{pb}(x_3)$ is weakly dominated by all points in $f(x_3, u_1, \mathcal{Y})_N$ and weakly dominates all points in $f(x_3, u_2, \mathcal{Y})_N$.} \label{scalarizations}
\end{figure}
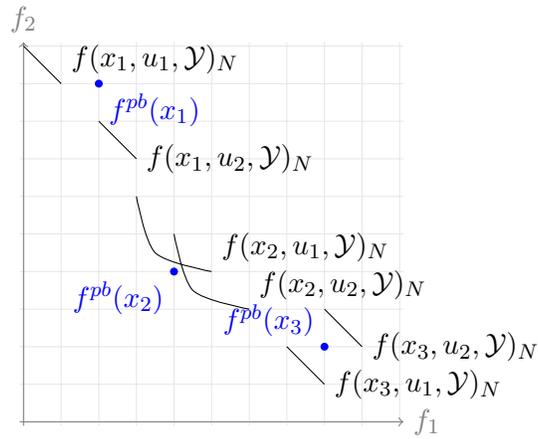

Figure \ref{scalarizations} shows that a componentwise interpretation can be too pessimistic as well as too optimistic. This is why this approach, that can be used to obtain conservative bounds for MRO problems, should not be applied to multicriteria \emph{adjustable} robust problems. In reality, only one scenario can occur and only one reaction to that is possible, in contrast to what this approach suggests.

Finally, we also look at the objective space image.
\begin{defi}
    We define the point-based min-max-min MARO-efficiency objective space image as
    \begin{align*}
    	im^f_{pb}:=\{f^{pb}(x) \mid x \text{ point-based min-max-min MARO-efficient}\}.
    \end{align*}
\end{defi}

\begin{lem} \label{pointfront}
	All points in $im^f_{pb}$ are nondominated.
\end{lem}

\begin{pf}
	This follows as all points in $im^f_{pb}$ are nondominated for a single stage multicriteria optimization problem (where the objective functions are multi stage functions). \qed
\end{pf}

\subsection{Comparison} \label{comparison}

\begin{table}
	\adjustbox{max width=1.3\textwidth, angle=90}{
			\renewcommand{\arraystretch}{1.5}
		\setlength{\tabcolsep}{10pt}
\begin{tabularx}{1.1\textheight}{>{\RaggedRight\arraybackslash}X>{\RaggedRight\arraybackslash}X>{\RaggedRight\arraybackslash}X>{\RaggedRight\arraybackslash}X}
	\toprule
	\textbf{computational concept} & \textbf{weighted sum MARO-efficiency} & \textbf{constraint MARO-efficiency} & \textbf{point-based min-max-min MARO-efficiency} \\
	\hline
	\hline
	\textbf{structure of objective space image} & Not all points in $im^f(\Lambda)$ have to be weakly nondominated (Remark~\ref{lambdafront}). & All points in $im_j^f(E)$ are weakly nondominated (Lemma~\ref{constraintfront}). & All points in $im^f_{pb}$ are nondominated (Lemma~\ref{pointfront}). \\
	\hline
	\textbf{bounds for all $u \in \mathcal{U}$} & $\sum_{i=1}^n \lambda_i f_i(x,u,y) \leq \mathcal{G}_\lambda$ (Remark~\ref{weightedbound}) & $f_i(x,u,y) \leq \varepsilon_i \; \forall i \neq j,$ $ f_j(x,u,y) \leq \mathcal{G}_{\varepsilon,j}$ (Remark~\ref{constraintbound}) & no better bounds than the trivial ones (Lemma~\ref{pointest}, Figure~\ref{badexample}) \\
	\hline
	\textbf{efficiency property for solutions which are strictly efficient with respect to the concept} & strictly multi-scenario MARO-efficient with respect to the preorder $\leqq^{\lambda, min}$ (Theorem~\ref{firstweightedlem}) & strictly multi-scenario MARO-efficient for the preorder $\leqq^{l}$ (Theorem~\ref{theostrictconst}) & - \\
	\hline
	\textbf{problem structure} & unchanged & changed & unchanged, but more variables ($n$ times) for middle and inner stage \\
	\bottomrule
\end{tabularx}
	}
	\caption{Comparison of computational concepts.}\label{Table} 
\end{table}

A comparison of the different methods from Section~\ref{approaches} is given in Table~\ref{Table}. The results for objective space image structure and bounds do also hold for weakly MARO-efficient solutions of the respective types, which can easily be computed for all compared MARO-efficiency types. It can be seen that weighted sum MARO-efficient solutions are only a valuable option if an explicit trade-off preference is known. In this case, there is no necessity to compute a Pareto front-like structure for decision aiding, and it is possible to profit from an unchanged structure of the problem. For general purposes, when no trade-off preference is known a priori, constraint MARO-efficiency seems to be the best option. It offers the possibility to compute kind of a Pareto front, to allow the use of common trade-off finding procedures, and gives the most interpretable bounds. However, the introduced constraints may change the problem structure and therefore lead to a harder problem. 
To visualize the differences, Figure~\ref{lambdanotconstr} shows an example of solutions which are weighted sum MARO-efficient but not constraint MARO-efficient and vice versa.
If solutions are strictly weighted sum or constraint MARO-efficient, the computed solutions can be understood as solutions in the sense of Section~\ref{Problem}.
\definecolor{mycolor}{RGB}{0,100,0}
\begin{figure}
	\parbox{0.5\textwidth}{
		\begin{tikzpicture}[scale=0.5, domain=0:10]
			\draw[very thin,color=lightgray] (-0.1,-0.1) grid (10.1,10.1);
			
			\draw[->,color=gray] (-0.1,0) -- (10.1,0) node[right] {$f_1$};
			\draw[->,color=gray] (0,-0.1) -- (0,10.1) node[above] {$f_2$};
			
			\draw plot coordinates {(2,7+3) (5,4+3) (10, 2+3)} node[above] {$f(x_1, u_1, \mathcal{Y})_N$};
			\draw plot coordinates {(2,6+3) (5,3+3) (10, 1+3)} node[below] {$f(x_1, u_2, \mathcal{Y})_N$};
			\draw[color=blue, opacity=0.8] plot[smooth] coordinates {(1, 5+3) (2,4+3) } node[above right] {$f(x_2, u_1, \mathcal{Y})_N$};
			\draw[color=blue, opacity=0.8] plot[smooth] coordinates {(1, 4.5+3) (2,3.5+3)} node[below right] {$f(x_2, u_2, \mathcal{Y})_N$};
			\draw[color=mycolor, opacity=0.8] plot coordinates {(8, 3) (9,2)} node[above left] {$f(x_3, u_1, \mathcal{Y})_N$};
			\draw[color=mycolor, opacity=0.8] plot coordinates {(8, 2) (9,1)} node[left] {$f(x_3, u_2, \mathcal{Y})_N$};
		\end{tikzpicture}
	}
	\begin{minipage}{0.5\textwidth}
	\begin{tikzpicture}[scale=0.5, domain=0:10]
		\draw[very thin,color=lightgray] (-0.1,-0.1) grid (10.1,10.1);
		
		\draw[->,color=gray] (-0.1,0) -- (10.1,0) node[right] {$f_1$};
		\draw[->,color=gray] (0,-0.1) -- (0,10.1) node[above] {$f_2$};
		
		\draw plot coordinates {(2,9) (6,3)} node[below left] {$f(x_1, u_1, \mathcal{Y})_N$};
		\draw plot coordinates {(2,5) (6,4)} node[above right] {$f(x_1, u_2, \mathcal{Y})_N$};
\draw[color=blue, opacity=0.8] plot coordinates {(2, 6) (6.5,3.2) } node[right] {\begin{tabular}{c} $f(x_2, u_1, \mathcal{Y})_N$  \\ $= f(x_2, u_2, \mathcal{Y})_N$ \end{tabular}\\ };
	\end{tikzpicture}
	\end{minipage}
	\caption{On the left side, $x_1$ is constraint MARO-efficient but not weighted sum MARO-efficient. On the right side, $x_1$ is weighted sum MARO-efficient but not constraint MARO-efficient.}\label{lambdanotconstr} 
\end{figure}
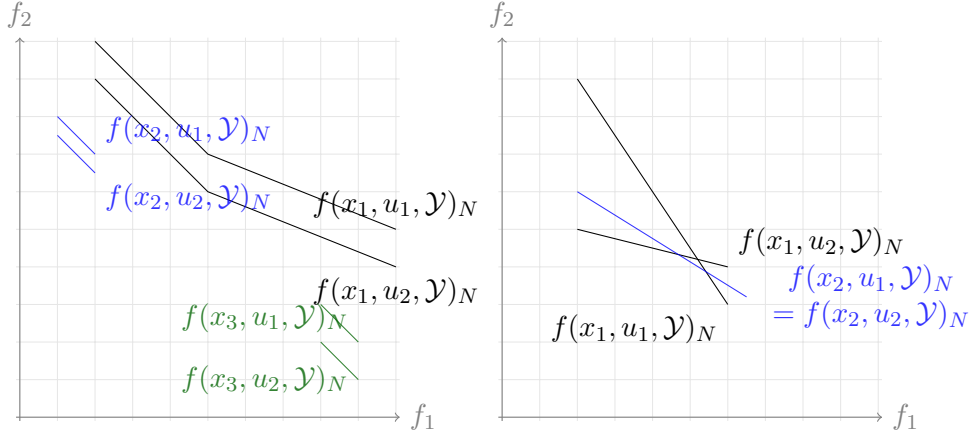

The last option is point-based min-max-min MARO-efficiency, which is not recommended, because of the lack of interpretability due to the non-existence of valuable bounds and the phenomenon that can be seen in Figure~\ref{scalarizations}. Of course, this cannot be compensated for by the objective space image structure and the unchanged, yet expanded, problem structure.

\section{Conclusion and Outlook} \label{conclusionandoutlook}
In this paper, we suggested different definitions of MARO-efficiency. Those from Section~\ref{Problem} are extensions of definitions for MRO problems. Although these extensions are reasonable and consistent with MCO under decision uncertainty and MRO, they are not immediately useful in practice, as it is not clear how to find such solutions. That is why we looked at alternative definitions in Section~\ref{approaches}, which are based on scalarization. We found that an $\varepsilon$-constraint inspired approach has the most desirable properties when compared to a weighted sum or a first-robustify approach. 

There are several directions for future research on this topic.
\begin{enumerate}
	\item The solution of general single criteria three stage problems. 
	\item Exploration of approaches other than the ones given in Section~\ref{approaches}, for example, a weighted Chebyshev scalarization inspired approach or an approach that is inspired by augmented $\varepsilon$-constraint method.
	\item Identification of conditions under which more relations between different efficiency definitions from Sections~\ref{Problem} and \ref{approaches} can be shown. 
	\item We assumed an understanding of $[\leqq/\leq/<]$ that is induced by the cone $\mathbb{R}_{[\geqq/\geq/>]}^n$. The results could be extended to other cones.
\end{enumerate}

\section*{Declarations of interest: }
none.

\section*{Acknowledgements}
This research did not receive any specific grant from funding agencies in the public, commercial, or not-for-profit sectors.

\section*{CRediT authorship contribution statement}
E. Halser: Conceptualization, Formal analysis, Methodology, Visualization, Writing - original draft.
E. Finhold: Writing - review \& editing.
N. Leith\"auser: Supervision.
J. Schwientek: Writing - review \& editing.
K. Teichert: Conceptualization, Writing - review \& editing.
K.-H. K\"ufer: Supervision.



\bibliographystyle{elsarticle-harv} 
\bibliography{all_sources}





\end{document}